\documentclass{elsart}

\usepackage{natbib}
\usepackage{amssymb,amsfonts, amsmath}
\usepackage[dvips]{graphicx}

\begin{document}

  \begin{frontmatter}

\title{The Space of Compatible Full Conditionals is a Unimodular Toric
  Variety}   

\author[abs]{Aleksandra B Slavkovic}
\author[ss]{and Seth Sullivant}
\address[abs]{Department of Statistics, Carnegie Mellon University}
\address[ss]{Department of Mathematics, University of California, Berkeley}

\begin{abstract}
The set of all $m$-tuples of 
compatible full conditional distributions on discrete random variables is an
algebraic set whose defining ideal is a unimodular toric ideal. We
identify the defining polynomials of these ideals with closed walks on a 
bipartite graph.  Our algebraic characterization provides a natural
generalization of the requirement that compatible conditionals have
identical odds ratios and holds regardless of the
patterns of zeros in the conditional arrays. 
\end{abstract}

\begin{keyword}
Odds ratios, Toric ideal, Unimodular
\end{keyword}

\end{frontmatter}



\section{Introduction}
Statisticians have long been interested in combining marginal
and conditional distributions in order to completely specify a joint
distribution.  Such models appear in spatial statistics
\citep{Besag1974}, analysis of contingency tables \citep{Bishop1975},
Bayesian prior elicitation 
contexts \citep{O'Hagan1998}, expert systems \citep{Cowell1999},
statistical disclosure
limitation \citep{Slavkovic2004b}, and generally in any area of applied
statistics where one
wishes to build global statistical models from local information about
subcollections of random variables.

Two fundamental theoretical questions which have been addressed in the
literature are the compatibility of conditionals and marginals, and
the uniqueness of the joint distribution when it exists. A collection 
of conditionals and marginals are compatible if there exists a joint
distribution with these conditionals and marginals.    The first major
result along these
lines is the Hammersley-Clifford theorem 
\citep{Besag1974} which establishes a connection between the positive joint
probability distribution and the full conditionals: it provides sufficient
conditions for compatibility of distributions in the setting of 
Markov-random fields. In the discrete case, the Hammersley-Clifford
theorem applies when there are
no cells with zero probability. Questions about the uniqueness of the joint
distribution given a compatible collection of marginals and
conditionals were addressed in \citet{Gelman1993} and the results
presented there were subsequently generalized in \citet{Arnold1999}.
\citet{Arnold1999} also address the 
question of whether the given set of densities are compatible; they
describe a variety of compatibility conditions for the case of
finitely many discrete random variables
and give algorithms for checking the compatibility of a set of
conditional probability densities. 


Let $X_1, \ldots, X_n$ be discrete random
variables where the set of possible states for each $X_i$ is in the
set of integers $[d_i] := \{1, 2, \ldots, d_i\}$.  The joint
probability distribution $p(x) = p(X_1 = i_1, \ldots, X_n = i_n)$
can be represented by a $d_1 \times \cdots \times d_n$ array of
nonnegative real numbers that sum to one.   
We define a \emph{full conditional} as a conditional multi-dimensional
array $p(x_A|x_B)$ where $A \cup B
= [n]$ and $A \cap B = \emptyset$ so that the conditional array
depends on all the random 
variables.  The precise statement of the compatibility problem for
full conditionals is the following:

\begin{prob}
Given partitions $(A_1, B_1),$ $\ldots, (A_m, B_m)$ of $[n]$ and
arrays $C^{1} \ldots C^{m}$ each of format $d_1 \times \cdots \times
d_n$, when is there 
a joint probability distribution $p(x)$ such that $p(x_{A_i} | 
x_{B_i}) = C^{i}$ for all $i$?
\end{prob}

The compatibility problem makes sense for any set of conditionals
and marginals, but we shall see a particularly simple combinatorial
solution when we restrict attention to compatibility among full
conditionals.

Given fixed values of the conditional and marginal arrays, the
problem of determining the existence of a 
joint probability distribution can be phrased as a feasibility
question for linear programs.  Since marginal and conditional arrays
are defined by linear and linear-fractional constraints in 
terms of the joint probabilities, the existence of a joint given the
conditionals can be decided by determining whether or not a
certain polytope is nonempty \citep{Arnold1999}.  While this
solution has many attractive features (the most important of which
is the ability to include extra parameters to account for
data that are ``almost compatible" \citep{Arnold1999}), it does not give an
intrinsic characterization of those sets of arrays which
correspond to compatible conditional distributions.  In
particular, using linear programming to determine whether or not a set
of conditionals  are
compatible requires the introduction of extra parameters.

In this paper, we develop intrinsic methods for
determining compatibility of conditional matrices, in the presence of
zero cell entries. We describe explicit algebraic restrictions on
the arrays $C^1, \ldots, C^m$ which guarantee that the
prescribed full conditional distributions are compatible.  The results
can be interpreted statistically in 
terms of conditions that replace generalized odds ratios to allow for
arbitrary patterns of zeros.  In this sense, our results provide a key
link between the statistical notions of odds ratios and tools from
algebraic geometry.

Our main result is the following theorem.
\begin{thm}\label{thm:main}
An $m$-tuple of $d_1 \times \cdots \times d_n$ arrays $C^1, \ldots,
C^m$ are  compatible full conditionals for the partitions $(A_1,B_1),
\ldots, (A_m,B_m)$ if and only if they satisfy the following four
conditions: 
\begin{enumerate}
\item  $C^i_{j_1 \cdots j_n} \geq 0$ for all $i$ and $j_1,  j_2, \ldots, j_n$.
\item  If $C^{i_0}_{j_1 \cdots j_n} = 0$ for some $i_0$, then
  $C^i_{j_1 \cdots j_n} = 0$ for all $i$. 
\item  For each $i$ the $B_i$ margin of $C^i$ is an array of all ones.
\item The entries of the arrays $C^i$ satisfy polynomial conditions
  which are in bijection with the induced circuits of a
  bipartite graph that depends on the $d_1, d_2, \ldots, d_n$ and
  $(A_1, B_1), \ldots, (A_m, B_m)$. 
\end{enumerate}
\end{thm}

The outline for the rest of the paper is as follows.  In the next
section, we provide a detailed description of our results in the case
of two random variables.  We give precise defining relations for the
space of compatible conditionals which generalize the condition that
odds ratios of consistent conditionals should be the same.  In
particular, we give necessary 
and sufficient polynomial conditions for two conditionals to be
compatible.  In the third section, we describe the theoretical tools
from  combinatorial commutative algebra
\citep{St} and discrete optimization \citep{sch} which we will need to
prove Theorem \ref{thm:main}.  In the fourth section 
we provide a general proof that the
space of compatible full conditionals is a unimodular toric variety.
The fifth section is devoted to a few
trivariate examples to  illustrate the main theorem.

\section{Two random variables}

In this section, we give a detailed explanation of our algebraic
results in the case of the bivariate compatible conditional problem.
We refer the reader to \citep{CLO} for the basics of algebraic geometry.
We relegate technical definitions and proofs of the main theorems to
the third and fourth sections.  
We denote the two potential
conditional arrays by
$$C_{ij} = p(X_1 = i | X_2= j) = \frac{p(X_1 = i, X_2 = j)}{p(X_2 =
  j)}\mbox{  and  }$$ $$ D_{ij} = p(X_2= j| X_1= i) = \frac{p(X_1 = i, X_2
  = j)}{p(X_1 = i)}.$$
This yields the following parametrization of the space of compatible
conditionals. 

\begin{prop}\label{thm:lawr}
Two nonnegative arrays $C$ and $D$ represent compatible conditionals for two
discrete random 
variables if and only if there are
nonnegative parameter arrays $P_{ij}$, $u_i$, $v_j$ such that
$$C_{ij} = P_{ij} v_j \mbox{  and  } D_{ij} = P_{ij}u_i,$$
\noindent $\sum_i C_{ij} = 1$ for all $j$ and $\sum_j D_{ij} = 1$ for
all $i$.
\end{prop}

Statistically, the array  $P_{ij}$ represent a joint
distribution which has $C$ and $D$ as its corresponding
conditionals.  Arrays $u_i$ and $v_j$ are marginals which when
combined with appropriate conditionals give $P_{ij}$. This
representation says that the set of all arrays $C$ and $D$ which are
compatible conditions is determined by allowing the arrays $P$, $u$
and $v$ to range over all nonnegative values which make $C$
and $D$ into conditional arrays.  The given
parametrization is by polynomial functions in the parameters, and
hence there exist polynomials 
in the $C_{ij}$ and $D_{ij}$ which vanish if and only if the pair $C$ and
$D$ belong to the (closure of the) space of compatible conditionals. 
In the literature of discrete optimization and toric ideals, the
parametrization defined in Proposition \ref{thm:lawr} is known as the
Lawrence lifting of the Segre variety.  The
equations which vanish on this parametrization are well understood.

\begin{thm}\citep[e.g.][Ch. 14]{St}
The defining equations for the Lawrence lifting of the Segre variety are in
bijection with circuits in the complete bipartite graph $K_{d_1,d_2}$.
Each circuit of  length $2r$ defines a binomial of degree $2r$ in the
$C_{ij}$ and $D_{ij}$. 
\end{thm}

This bijection is described as follows.  Let $(i_1,j_1, i_2, j_2, \ldots,
i_r,j_r,i_1)$ be a circuit in $K_{d_1,d_2}$ of length $2r$.  This
circuit produces a binomial of degree $2r$:
$$ C_{i_1j_1} D_{i_2 j_1} C_{i_2 j_2} \cdots D_{i_1 j_r} -  D_{i_1j_1}
C_{i_2 j_1} D_{i_2 j_2} \cdots C_{i_1 j_r}.$$
Note that the degree four relations which are produced by this
construction are precisely the relations which say that all the odds
ratios of the matrices $C$ and $D$ must be the same.  For example, the
circuit of length 4 in $K_{2,2}$ produces the binomial: 
$$C_{11}C_{22}D_{12}D_{21} - C_{12}C_{21}D_{11}D_{22}.$$
For positive $C$ and $D$, this binomial is zero if and only if the
following odds ratios are equivalent:
$$\frac{C_{11}C_{22}}{C_{12}C_{21}} = \frac{D_{11}D_{22}}{D_{12}D_{21}}.$$

It is well known that for positive matrices $C$ and $D$, compatibility
of conditional matrices is  
equivalent to equality of the odds ratios \citep{Arnold1999}. In
algebro-geometric terms,  this says that the degree four
polynomials produced by the circuit construction define
the toric variety of interest set theoretically in the strictly
positive orthant. 
However, this statement becomes false when some of the entries in $C$
and $D$ are allowed to be zero.   These polynomials, in
addition to addressing the compatiblity of arrays $C$ and
$D$ in the presence of zero entries, 
also should be useful for determing the similarity of linear contrasts
in the $C$ and $D$ matrices when there are zero cells
\citep{Bishop1975}.  Thus, these polynomials could prove useful for
the analysis of incomplete contingency tables.   

\begin{exmp}
Consider the matrices
$$ C = \left( \begin{array}{ccc}
\frac{1}{2} & \frac{1}{2} & 0 \\
0 & \frac{1}{2} & \frac{1}{2} \\
\frac{1}{2} & 0 & \frac{1}{2} \end{array} \right) \mbox{   and   }
D = \left( \begin{array}{ccc}
\frac{1}{3} & \frac{2}{3} & 0 \\
0 & \frac{1}{3} & \frac{2}{3} \\
\frac{1}{3} & 0 & \frac{2}{3} \end{array} \right). $$
All of the odds ratios of these matrices are the same:  they
are either simultaneously zero or undefined.  In other words, they
saitisfy all 9  polynomials which are described by circuits in
$K_{3,3}$ of length 4.   However, the matrices do
not represent compatible conditionals since the binomial given by the length 6
circuit $(1,1,2,2,3,3,1)$: 
$$ C_{11}C_{22}C_{33}D_{12}D_{23}D_{31} -
C_{12}C_{23}C_{31}D_{11}D_{22}D_{33}$$ 
is equal to $1/108$ and not zero.  These incompatible conditionals appear
in \cite{Arnold1999}.
\end{exmp}

These polynomials give a parameter free
method for determining if two arbitrary conditional
arrays are compatible.  

To make the connection to the next section more lucid, we will
introduce an alternate graph theoretic representation of the circuit
binomials described above.  This more complicated description
is the one that generalizes for arbitrary sets  of full conditionals
in the next section. 

\begin{defn}
We define the bipartite graph $G_{d_1,d_2}$ to be the graph with
$d_1d_2 + d_1 + d_2$ vertices as follows. 
\begin{enumerate}
\item  There are $d_1 d_2$ vertices labeled $v_{ij}$ where $i \in [d_1]$ and $j \in [d_2]$.
\item  There are $d_1$ vertices labeled $w_k$ where $k \in [d_1]$.
\item  There are $d_2$ vertices labeled $u_l$ where $l \in [d_2]$.
\item  There is an edge between $v_{ij}$ and $w_k$ if and only if $i = k$.
\item There is an edge between $v_{ij}$ and $u_l$ if and only if $j = l$.
\item  There are no other edges in $G_{d_1, d_2}$.
\end{enumerate}
\end{defn}

\begin{figure}\label{g23}
\begin{center}\includegraphics{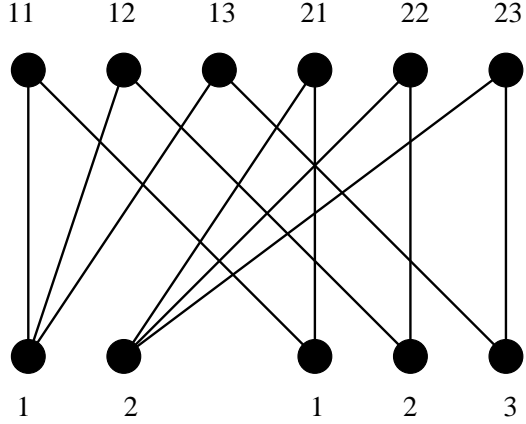}
\caption{The bipartite graph $G_{2,3}$}
\end{center}
\end{figure}

The graph $G_{2,3}$ is pictured in Figure 1.  Associated to any
circuit in $G_{d_1,d_2}$ of length $2r$ we get a 
binomial of degree $r$ as follows.  After permuting indices, circuits
in $G_{d_1,d_2}$ have the form 
$$( i_1j_1,i_1,i_1j_2,j_2,i_2j_2, \ldots,  j_r,i_rj_r)$$
\noindent where $i_1 = i_r$ and $j_1 = j_r$.  From this circuit we
recover the binomial 
$$ C_{i_1j_1} D_{i_2 j_1} C_{i_2 j_2} \cdots D_{i_1 j_r} -  D_{i_1j_1}
C_{i_2 j_1} 
D_{i_2 j_2} \cdots C_{i_1 j_r},$$

which was the same binomial from our previous description.
Putting all these ideas together we have the following special case of
Theorem \ref{thm:main}.

\begin{cor}
A pair of matrices $C$ and $D$ are compatible conditionals if and only
if they satisfy the following properties:
\begin{enumerate}
\item $C \geq 0 $ and $D \geq 0$,
\item for all $i$ and $j$, $C_{ij} = 0$ if only if $D_{ij} = 0$,
\item $\sum_i C_{ij} = 1$ for all $j$ and $\sum_j D_{ij} = 1$ for all $i$,
  and
\item $C$ and $D$ satisfy the circuit polynomials associated to the
  circuits in the graph $G_{d_1,d_2}$.
\end{enumerate}
\end{cor}

\section{Unimodular toric varieties}

In this section, we will describe the mathematical objects that we
will need to describe the space of compatible conditionals.  First, we
show that the space of compatible conditionals is a parametrized
algebraic set.  The parametrization is given by monomials and this
leads naturally to the definition of \emph{toric varieties} and
\emph{toric ideals}.  Then we introduce the notions of graphical
unimodular matrices and unimodular toric varieties, which are the
natural algebraic objects for representing the space of compatible
conditionals. 

A key preliminary observation is the following:

\begin{thm}
A set of full conditional matrices $C^1, \ldots, C^m$ are compatible
conditionals corresponding to the partitions $(A_1, B_1), \ldots,
(A_m,B_m)$ 
if and only if there are matrices of probabilities $U^1, \ldots, U^m$ of
the appropriate sizes and a $d_1 \times \cdots \times d_n$ matrix
of probabilities $P$ such that 
$$P_{j_1, \ldots, j_n} = C^i_{j_1,\ldots, j_n} \cdot
 U^i_{j_{k_1},\ldots, j_{k_s}}$$
for all $i$ where $\{k_1 < \ldots < k_s\} = B_i.$
\end{thm}

\begin{pf}
This is the definition of
compatibility together with Bayes' rule.  The matrices $U^i$ are then
the ``missing''$B_i$-marginal distributions. \qed
\end{pf}

Since we are interested in intrinsic characterizations of the
compatible conditionals, it makes sense to rearrange this expression
to deduce an equation for the $C^i$ in terms of the missing joint
and conditional distributions $P$ and $U^i$.  Thus we deduce the
following corollary. 

\begin{cor}\label{cor:tor}
The set of all $m$-tuples of consistent full conditional arrays $C^1$,
$\ldots$, $C^m$ is 
described parameterically as
$$C^i_{j_1 \cdots j_n}  = P_{j_1  \cdots j_n} U^i_{j_{k_1} \cdots j_{k_s}}.$$
where the multiplication is componentwise.
\end{cor}

At first glance, the space of conditionals we have described in
Corollary \ref{cor:tor} looks highly over-parameterized:
shouldn't the $U$ parameters be related to the $P$ parameters by
marginalization?

We claim that we can first allow the parameters $P$ and $U^i$ to be
arbitrary real arrays, consider the variety which is defined by this
parameterization, and then intersect this variety with the appropriate
product of probability simplices.  The result of this procedure will
be the same set of conditional arrays which is obtained by first
restricting the parameters so that $U^i = p(x_B)$ and then considering
the resulting parametrically described variety.  To see this, observe
that the only way that $C^i$ could be a 
conditional array is if $U^i$ is the corresponding marginal $p(x_B)$,
which again follows by Bayes' rule. 
 
Let $p(x_{A_1} | x_{B_1})$ and $p(x_{A_2} | x_{B_2})$ be two full
conditional matrices.  Note that if $B_2 \subseteq B_1$, then they
are compatible if and only if $p(x_{A_1} | x_{B_1})$ is obtained
from $p(x_{A_2} | x_{B_2})$ by conditioning on the variables $A_2
\setminus A_1$.  Hence when determining the space of consistent
conditionals for $p(x_{A_1} | x_{B_1})$, $p(x_{A_2} | x_{B_2})$,
$\ldots$, $p(x_{A_m} | x_{B_m})$, we may assume that there are no
containment relations among the $A_i$ or $B_i$.  That is, the
collection $\Delta = \{B_1, \ldots, B_m \}$ can be taken to be the
facets of a simplicial complex.

We will work in the polynomial ring
$$R = \mathbb{R}[C^i_{j_1, \ldots,
j_n}| 1 \leq  i \leq m \mbox{ and } 1 \leq j_s \leq d_s, \mbox{ for
    all}  s]$$
 which is a polynomial ring in $\mathbf{D}=md_1d_2 \cdots d_n$
indeterminates.  Each indeterminate in the ring $R$ corresponds to an
entry in one of the $m$ conditional matrices.  By Corollary
\ref{cor:tor}, the algebraic variety in $\mathbb{R}^\mathbf{D}$,
which, when intersected with a product of probability simplices,
yields the space of consistent conditional, is given by the monomial
parameterization 
$$ \label{eq:param} C^i_{j_1, \ldots, j_n} = P_{j_1, \ldots, j_n}
U^i_{j_{k_1}, \ldots, j_{k_s}}$$ 
where $\{ k_1 < k_2 \cdots < k_s \} = B_i$ and the $P$'s and $U$'s are
allowed to range over all of the real numbers. 

\begin{defn}
An algebraic variety which is given by a monomial parameterization is
called a \emph{toric variety}. 
\end{defn} 

To describe the polynomials that define any toric variety amounts to
understanding which products of the parameterizing monomials are equal
to which other products of the parameterizing monomials.  That is, the
\emph{toric ideal} which defines the toric variety is generated by
binomials (differences of monomials).  

\begin{defn}
The toric ideal which is the vanishing ideal of the space of
compatible conditionals is denoted $I_\Delta$ which  we call the
\emph{compatibility ideal}.  
\end{defn}

Furthermore, understanding which monomials are equal to each other
can be studied by understanding the integer kernel of an associated
matrix.  A standard reference for the theory of toric ideals is \citep{St}.

\begin{defn}\label{def:tor}
Let $A \in  \mathbb{Z}^{d \times n}$ be an integer matrix.  The toric
variety associated to $A$ is defined by the parameterization 
$$ x_i = t_1^{a_{1i}} t_2^{a_{2i}} \cdots t_d^{a_{di}}.$$
\noindent  The toric ideal $I_A$ which defines this parameterized
variety is generated by the (infinitely many) binomials 
$$ I_A = \left<  x^u - x^v | Au = Av \right>. $$
\end{defn}

The presentation of a toric ideal $I_A$ in Definition \ref{def:tor} is
given with infinitely many ideal generators.  However, the Hilbert basis
theorem implies that only finitely many of these polynomials 
are actually needed in a generating set of the ideal.  In
many situations, the difficulty of studying toric ideals amounts to
finding good combinatorial descriptions of this finite set of
generators.  In the case of the compatibility ideals $I_{\Delta}$
there is a simple combinatorial structure called 
\emph{unimodularity} which can be used to describe the generators.

\begin{defn}
A toric ideal $I_A$ is called \emph{unimodular} if every reduced Gr\"obner
basis of $I_A$ consists of squarefree binomials (i.e. there are no
squared variables in the leading terms of each
binomial in the reduced Gr\"obner basis).  Equivalently, $I_A$ is
unimodular if all the non-zero minors of $A$ have the same absolute value.
Such a matrix $A$ is also called unimodular. 
\end{defn}

The following is an important special family of unimodular matrices
(see, for example \citet{sch}).

\begin{prop}
Let $A$ be a matrix with the following properties:
\begin{enumerate}
\item  All the entries of $A$ are either zero or one,
\item  there are precisely two 1's in each column of $A$, and
\item  there is a partition of the rows of $A$ into two sets $U$ and
  $V$ such that each column of $A$ has exactly one nonzero entry with
  row index $u \in U$ and one nonzero entry with row index $v \in V$.
\end{enumerate}
Then $A$ is a unimodular matrix and is called a \emph{graphical
  unimodular matrix}.
\end{prop}

The reason for the name ``graphical'' is due to the fact that any
graphical unimodular matrix $A$ is the vertex-edge incidence matrix of a
bipartite graph $G_A$.  The parition of the vertices of $G_A$
corresponds to the partition of the rows of $A$.  Notice that vectors
$v \in \ker(A)$ correspond to the union of cycles in the graph $G_A$.

\begin{defn}
Let $G$ be any graph.  An \emph{induced circuit} of $G$ is a circuit
of $G$ which does not have a chord (i.e. a ``short cut'') in $G$.
\end{defn}

In the complete bipartite graph $K_{3,3}$ all of the cycles of length
4 are induced, 
whereas none of the six cycles are induced because there is a chord
cutting across them.  Toric ideals which are
presented by a graphical unimodular matrix have a simple combinatorial
description for their minimal generators in terms of the induced
circuits of the associated bipartite graph.  A version of the
following result can be found in \citep{AT}.

\begin{prop}
Let $A$ be a graphical unimodular matrix, and $G_A$ the associated
bipartite graph.  Let the $U$ and $V$ be the partition of the vertices
of $G_A$ where $U = \{u_1, \ldots, u_{l_1}\}$ and $V = \{v_1, \ldots,
v_{l_2}\}$.  For each circuit $c = (u_{i_1},v_{j_1}, u_{i_2}, v_{j_2},
\ldots, u_{i_r}, v_{j_r}, u_{i_1})$ we associate the circuit binomial
$$f_c = x_{u_{i_1}v_{j_1}} x_{u_{i_2}v_{j_2}} \cdots
x_{u_{i_r}v_{j_r}} -
x_{u_{i_2}v_{j_1}} x_{u_{i_3}v_{j_2}} \cdots x_{u_{i_1}v_{j_r}}.$$
Then $I_A$ is minimally generated by the circuit binomials
corresponding to the induced circuits of $G_A$.  That is,
$$ I_A = \left< f_c \mid c \mbox{ is an induced circuit of } G_A \right>. $$
\end{prop}

\section{The main algebraic result}

Now that we have reviewed all the mathematical facts we need about
unimodular toric ideals, we are in a position to prove that the
compatibility ideal is a unimodular toric ideal.  To do this, we
construct the matrix that represents this toric ideal.

Let $d = (d_1, \ldots, d_n)$ be the integer vector corresponding to
the dimensions of the tables let $\Delta = \{B_1, \ldots, B_m \}$ be
the simplicial complex which represents the compatible conditionals.
For each facet $B_i$ of $\Delta$, denote by $d_{B_i}$ the vector 
$$d_{B_i} = (d_{k_1} , \ldots, d_{k_s})$$ 
where $B_i = \{k_1, \ldots, k_s\}$ and denote by $D_{B_i}$ the product
$$D_{B_i} = \prod_{k_i \in F_i} d_{k_i}.$$

We denote by $A_{\Delta,d}$ the matrix that represents the
compatibility ideal corresponding to this data.  This matrix is a
$$ (d_1 d_2 \cdots d_n + \sum_{i = 1}^m d_{B_i}) \times m d_1 d_2
\cdots d_n $$
matrix.  The rows are naturally grouped into $m+1$ blocks
and the columns are natrually grouped into $m$ blocks.  Each column is
labelled by an integer $i \in [m]$ and an $n$-tuple $(j_1, \ldots,
j_n)  \in [d_1] \times [d_2] \times \cdots \times [d_n]$.  The first
block of rows is labelled by the integer $0$ and $n$-tuple $(j_1, \ldots,
j_n)  \in [d_1] \times [d_2] \times \cdots \times [d_n]$.  Each of the
remaining blocks of rows is labelled by and integer $i$ and an
$|B_i|$-tuple $(j_{k_1} , \ldots, j_{k_s}) \in [d_{k_1}] \times \cdots
\times [d_{k_s}]$ where $B_i = \{k_1, \ldots, k_s \}$.  

The entries of this matrix are all ones and zeros.  There is a one in
a particular entry with row indexed by the data $i^1$, $(j^1_1,
\ldots,j^1_d)$ and column indexed by the data $i^2$, $(j^2_1,
\ldots,j^2_d)$ if and only if it satisfies the following rules:

\begin{enumerate}
\item if the row label is $0$ and $ (j^1_1, \ldots,j^1_d) =
  (j^2_1,\ldots,j^2_d)$ or
\item if the row label is $i^1 > 0$, $i^1 = i^2$, and $ (j^1_{k_1},
\ldots,j^1_{k_s}) =  (j^2_{k_1},\ldots,j^2_{k_s}) $ where $B_{i^1} =
(k_1, \ldots, k_s).$
\end{enumerate}

\begin{thm}\label{thm:circuit}
The matrix $A_{\Delta,d}$ is a graphical unimodular matrix and it
represents the toric variety of the space of compatible conditionals.  Hence,
the compatibility ideal is generated by the induced circuit binomials
in the associated bipartite graph.
\end{thm}

\begin{pf}
To see that $A_{\Delta,d}$ represents the toric variety of compatible
conditionals amounts to identifying the rows of the matrix with a
parameter, and the columns of the matrix with an entry in the
conditional matrix.  The labelling for the columns $i, (j_1, \ldots,
j_n)$ naturally corresponds to the indeterminate $C^i_{j_1, \ldots,
  j_n}$.  A row in the first block, with $i = 0$, corresponds to the
$P_{j_1, \ldots, j_n}$ parameters.  A row in any of the blocks with $i
>0$ corresponds to the parameter $U^i_{j_{k_1}, \ldots, j_{k_s}}.$
Note that be the description of the matrix $A_{\Delta, d}$, in the
column corresponding to $C^i_{j_1,\ldots, j_n}$ there are precisely
two nonzero entries: one in the first block which corresponds to
$P_{j_1, \ldots, j_n}$ and one in block $i$, which corresponds to
$U^i_{j_{k_1}, \ldots, j_{k_s}}$.  Thus, this represents the
parameterization
$$C^i_{j_1 \cdots j_n}  = P_{j_1  \cdots j_n} U^i_{j_{k_1} \cdots
  j_{k_s}}$$
as desired.

Now we wish to show the unimodularity of $A_{\Delta, d}$, however,
this is an immediate consequence of the preceding argument.  Each
column has precisely two ones: one with index $i = 0$
and one with index $i> 0$.  Hence, the value of $i$ determines the
partition of the rows to deduce the structure of the underlying
bipartite graph. \qed
\end{pf}

Now we will explicitly describe the graph encoded by the
matrix $A_{\Delta,d}$.

\begin{defn}
The graph $G_{\Delta,d}$ associated to the compatibility ideal has
$$ (d_1 d_2 \cdots d_n + \sum_{i = 1}^m d_{F_i})$$
\noindent vertices and $ m d_1 d_2 \cdots d_n$ edges:  these are the
number of rows and columns of $A_{\Delta,d}$ respectively.  The
vertices are partitioned into two class: those labelled with $i = 0$
and some $(j_1, \ldots, j_n)$, and those labelled with $i>0$ and some
$(j_{k_1}, \ldots, j_{k_s})$.  Vertices with $i= 0$ are incident only
to vertices with $i >0 $ and conversely.  In particular, the vertex $i
=0$, $(j^1_1, \ldots, j^1_n)$ is incident to $i>0$, $(j^2_{k_1},
\ldots, j^2_{k_s})$ if and only if $(j^1_{k_1} \ldots, j^1_{k_s}) =
(j^2_{k_1}, \ldots, j^2_{k_s})$.
\end{defn}

To close this section, we now provide the proof of Theorem
\ref{thm:main}.

{\bf PROOF of Theorem \ref{thm:main}:} 
Conditions (1)-(3) of Theorem \ref{thm:main} are clearly necessary.
Condition (4) is simply the expression of Theorem \ref{thm:circuit}
above.  What remains to show is that this conditions are, in fact,
sufficient.  Suppose that a particular realization of conditional
arrays $C^1, \ldots, C^m$ satisfy conditions (1) - (4).  We wish to
show that there exists parameter matrices $P$ and $U^1, \ldots, U^m$,
which have $C^1, \ldots, C^m$ as their image.  By conditions (1), (3),
and (4) together with a result of
\citet{GMS}, this can happen if and only if the columns of $A_{\Delta,d}$
which are indexed by the support of $C^1, \ldots, C^m$ form a
\emph{nice facial subset} of the columns of $A_{\Delta,d}$.  The nice
facial subsets of 
the columns of $A_{\Delta, d}$ are precisely those sets of columns
which are obtained by taking a collection of rows $R$ of $A_{\Delta,d}$,
and taking exactly those columns of $A_{\Delta,d}$ which have zeros in
the rows $R$.  However, the support of $C^1, \ldots, C^m$ is a nice
facial subset of the columns of $A_{\Delta,d}$ because of condition
(2):  the rows of $A_{\Delta,d}$ to choose are precisely those rows
labelled by $(0, (j_1, \ldots, j_n))$ for all the indices $(j_1,
\ldots, j_n)$ where the $C^i_{j_1, \ldots, j_n}$ are collectively
zero. \qed

\section{Trivariate examples}

We will conclude our paper with some examples on three binary random
variables to illustrate the construction and the types of polynomials
that appear.

\begin{exmp}
Consider the compatibility ideal associated with three binary random
variables and suppose
that $\Delta = \{\{1\},\{2\},\{3\}\}$.  In other words, we are
considering the compatibility of $p(x_1,x_2|x_3), p_(x_1,x_3|x_2)$ and
$p(x_2,x_3|x_1)$.  The matrix, $A_{\Delta,(2,2,2)}$ has 14 rows and 24
columns.  It is the zero/one matrix

\tiny
$$\left(
\begin{array}{cccccccc|cccccccc|cccccccc}
1 & 0 & 0 & 0 & 0 & 0 & 0 & 0 & 1 & 0 & 0 & 0 & 0 & 0 & 0 & 0 &
1 & 0 & 0 & 0 & 0 & 0 & 0 & 0 \\
0 & 1 & 0 & 0 & 0 & 0 & 0 & 0 & 0 & 1 & 0 & 0 & 0 & 0 & 0 & 0 & 0 &
1 & 0 & 0 & 0 & 0 & 0 & 0  \\
0 & 0 & 1 & 0 & 0 & 0 & 0 & 0 & 0 & 0 & 1 & 0 & 0 & 0 & 0 & 0 & 0 & 0 &
1 & 0 & 0 & 0 & 0 & 0 \\
0 & 0 & 0 & 1 & 0 & 0 & 0 & 0 & 0 & 0 & 0 & 1 & 0 & 0 & 0 & 0 & 0 & 0 & 0 &
1 & 0 & 0 & 0 & 0  \\
0 & 0 & 0 & 0 & 1 & 0 & 0 & 0 & 0 & 0 & 0 & 0 & 1 & 0 & 0 & 0 & 0 & 0
& 0 & 0 &1 & 0 & 0 & 0 \\ 
0 & 0 & 0 & 0 & 0 & 1 & 0 & 0 & 0 & 0 & 0 & 0 & 0 & 1 & 0 & 0 & 0 & 0
& 0 & 0 & 0 & 1 & 0 & 0 \\
0 & 0 & 0 & 0 & 0 & 0 & 1 & 0 & 0 & 0 & 0 & 0 & 0 & 0 & 1 & 0 & 0 & 0
& 0 & 0 & 0 & 0 & 1 & 0 \\
0 & 0 & 0 & 0 & 0 & 0 & 0 & 1 & 0 & 0 & 0 & 0 & 0 & 0 & 0 & 1 & 0 & 0
& 0 & 0 & 0 & 0 & 0 & 1 \\
\hline
1 & 0 & 1 & 0 & 1 & 0 & 1 & 0 & 0 & 0 & 0 & 0 & 0 & 0 & 0 & 0  & 0 & 0
& 0 & 0 & 0 & 0 & 0 & 0 \\
0 & 1 & 0 & 1 & 0 & 1 & 0 & 1  & 0 & 0 & 0 & 0 & 0 & 0 & 0 & 0  & 0 &
0 & 0 & 0 & 0 & 0 & 0 & 0 \\
\hline
0 & 0 & 0 & 0 & 0 & 0 & 0 & 0 &  1 & 1 & 0 & 0 & 1 & 1 & 0 & 0  & 0 & 0 &
0 & 0 & 0 & 0 & 0 & 0 \\ 
0 & 0 & 0 & 0 & 0 & 0 & 0 & 0 & 0 & 0 & 1 & 1 & 0 & 0 & 1 & 1  & 0 & 0
& 0 & 0 & 0 & 0 & 0 & 0 \\
\hline
0  & 0 & 0 & 0 & 0 & 0 & 0 & 0 & 0  & 0 & 0 & 0 & 0 & 0 & 0 & 0 & 1 &
1 & 1 & 1 & 0 & 0 & 0 & 0 \\
0 & 0 & 0 & 0 & 0 & 0 & 0 & 0 & 0  & 0 & 0 & 0 & 0 & 0 & 0 & 0 & 0 & 0
& 0 & 0 & 1 & 1 & 1 & 1 \\
\end{array} \right). $$
\normalsize
The graph corresponding to this compatibility problem has 14 vertices:
8 of one type and 6 of the other, and 24 edges between them.  It is
pictured in Figure 2.

\begin{figure}[h!]\label{fig:gr}
\begin{center} \includegraphics{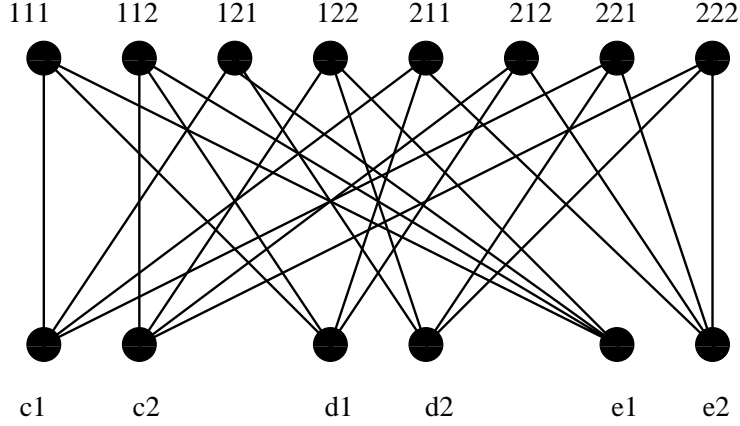}
\caption{The compatibility graph for $p(x_1,x_2|x_3), p(x_1,x_3|x_2), p(x_2,x_3|x_1)$}
\end{center}
\end{figure} 

There are three types of induced circuits in this graph, yielding three
types of binomial generators of the compatibility ideal $I_\Delta$.
We use the indeterminates $C_{ijk}$, $D_{ijk}$ and $E_{ijk}$ to denote
the corresponding entries in the conditional arrays.
For instance, the circuit $(111,c_1,211,d_1,111)$ in the graph is
induced and from it we deduce that the quadratic binomial
$$C_{111}D_{211} - D_{111}C_{211}$$
is a minimal generator of the compatibility ideal $I_{\Delta}$.
The induced circuit \\ $(111,c_1,221,e_2, 212, d_1,111)$ produces the
binomial
$$C_{111}E_{221}D_{212} - C_{221}E_{212}D_{111}$$
which is also a minimal generator of $I_\Delta$.  In total,
there are 12 induced circuits of length 4, 8 induced circuits of
length 6, and  60 induced circuits of length 8.  Hence $I_\Delta$ is
minimally generated by 12 quadric, 8 cubics, and 60 quartics.
\end{exmp}




\begin{exmp}
Consider the compatibility associated to 3 binary random variables
with the simplicial complex $\Delta = \{\{1,2\},\{1,3\}, \{2,3\}\}$.
That is, we are considering the compatibility of $p(x|y,z)$,
$p(y|x,z)$, and $p(z|x,y)$.  Use the indeterminates $C_{ijk}$,
$D_{ijk}$, and $E_{ijk}$ to represent entries in the corresponding
conditional arrays.  The ideal $I_{\Delta}$ is generated
by 28 binomials which fall into 4 equivalence classes modulo the
symmetry of the cube.  In this case, the binomial minimal generators
of $I_{\Delta}$ are in bijection with the circuit in the edge graph of
the cube.  The four different types of circuit on the
cube are indicated in figure 3.  

\begin{figure}[h!]\label{fig:cube}
\begin{center} \includegraphics{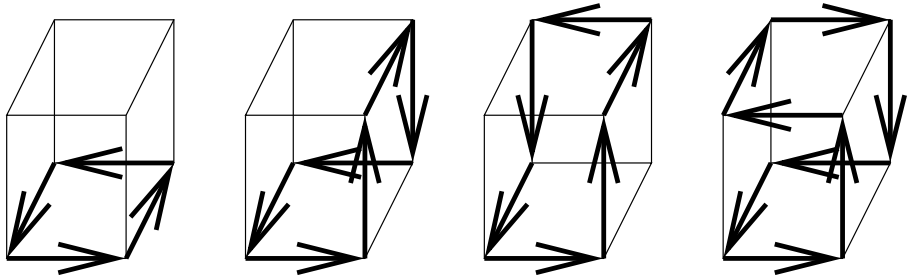}
\caption{The combinatorial types of circuits on the cube}
\end{center}
\end{figure}

These four types
of circuits yield four different symmetry classes of binomials in the
ideal $I_\Delta$ which are represented by the
following binomials:

$$
C_{111}C_{221}D_{121}D_{211} - 
C_{121}C_{211}D_{111}D_{221}
$$
$$
C_{111}D_{211}E_{221}D_{222}C_{212}E_{112} -
C_{211}D_{221}E_{222}D_{212}C_{112}E_{111},
$$
$$
C_{111}D_{211}E_{221}C_{222}D_{122}E_{112} -
C_{211}D_{221}E_{222}C_{122}D_{112}E_{111},
$$
$$
C_{111}D_{211}C_{221}E_{121}C_{122}D_{222}C_{212}E_{112} -
C_{211}D_{221}C_{121}E_{122}C_{222}D_{212}C_{112}E_{111}.
$$

Note that the first polynomial in this list has a natural statistical
interpretation:  it corresponds to the
equality of the odds ratios for $p(x|y,z)$ and $p(y|x,z)$.
\end{exmp}

While there are statistical
interpretations for some of the circuit polynomials which arise,
understanding these 
polynomial expressions in relation to known statistical quantities remains
an open problem.  In particular, determining how tools
from the
analysis of contingency tables such as generalized odds ratios and
linear contrast relate to the polynomial constraints we have derived
could prove useful for inference on incomplete contingency tables.

\section*{Acknowledgments}
Aleksandra Slavkovic was supported in part by National
Science Foundation Grant No. EIA-0131884 to the National Institute of
Statistical Sciences. 

\bibliographystyle{elsart-harv}

\end{document}